\magnification=1200 \baselineskip=12.3pt

\def\absdeter#1{\displaystyle\left\|\matrix{#1}\right\|}
\def\and{\quad\hbox{and}\quad}
\def\bC{{\bf C}}
\def\bQ{{\bf Q}}
\def\bR{{\bf R}}
\def\bx{{\bf x}}

\def\bZ{{\bf Z}}
\def\cC{{\cal{C}}}
\def\cE{{\cal{E}}}
\def\cO{{\cal{O}}}
\def\cqfd{\unskip\kern 6pt\penalty 500
  \raise -2pt\hbox{\vrule\vbox to 10pt{\hrule width 4pt
  \vfill\hrule}\vrule}\par}
\def\deter#1{\displaystyle\left|\matrix{#1}\right|}
\def\dZ#1{\left\{#1\right\}}
\def\GL{{\rm GL}}
\def\proof{\medskip\noindent{\it Proof.}\ }
\def\section#1{\bigskip\noindent{\bf #1} \medskip}
\def\td{\tilde{d}}
\def\tuy{\tilde{\uy}}
\def\trace{{\rm trace}}
\def\uT{{\bf T}}
\def\uw{{\bf w}}
\def\ux{{\bf x}}
\def\uy{{\bf y}}
\def\uz{{\bf z}}


\def\ADQZ{1}
\def\BT{2}
\def\Ca{3}
\def\DSa{4}
\def\DSb{5}
\def\LR{6}
\def\Lo{7}
\def\Ro{8}
\def\Sc{9}
\def\Wi{10}

\centerline{\bf Approximation to real numbers by cubic algebraic
integers I}

\bigskip
\centerline{ by Damien ROY}

\footnote{}{\sevenbf 2000 Mathematics Subject Classification:
\sevenrm Primary 11J04; Secondary 11J13, 11J82 }

\footnote{}{\sevenrm Work partly supported by NSERC and CICMA.}


\section{1. Introduction}

The study of approximation to a real number by algebraic numbers
of bounded degree started with a paper of E.~Wirsing [\Wi] in
1961. Motivated by this, H.~Davenport and W.~M.~Schmidt considered
in [\DSb] the analogous inhomogeneous problem of approximation to
a real number by algebraic integers of bounded degree.  They
proved a result that is optimal for degree $2$ and a general
result which is valid for any degree. We shall be concerned here
with the case of degree $3$, the first case for which the optimal
exponent of approximation is not known.

To state their results, we define the height $H(P)$ of a
polynomial $P\in \bR[T]$ to be the maximum of the absolute values
of its coefficients, and define the height $H(\alpha)$ of an
algebraic number $\alpha$ to be the height of its irreducible
polynomial over $\bZ$.  We also denote by $\gamma=(1+\sqrt{5})/2$
the golden ratio.

Davenport and Schmidt showed that, for any real number $\xi$ which
is neither rational nor quadratic irrational, there exist a
positive constant $c$ depending only on $\xi$ and arbitrarily
large real numbers $X$ such that the inequalities
$$
|x_0| \le X,
  \quad
|x_0\xi-x_1| \le cX^{-1/\gamma},
  \quad
|x_0\xi^2-x_2| \le cX^{-1/\gamma} \eqno(1.1)
$$
have no solutions in integers $x_0$, $x_1$, $x_2$, not all zero
(see Theorem 1a of [\DSb]).  By an argument of duality based on
geometry of numbers, they deduced from this that, for such a real
number $\xi$, there exist another constant $c>0$ and infinitely
many algebraic integers $\alpha$ of degree at most $3$ which
satisfy
$$
0< |\xi-\alpha| \le  c H(\alpha)^{-\gamma^2} \eqno(1.2)
$$
(see Theorem 1 of [\DSb]).

At first, it would be natural to expect that the first statement
holds with any exponent $>1/2$ in place of $1/\gamma \simeq 0.618$
in the inequalities (1.1), since this ensures that the volume of
the convex body defined by these inequalities tends to zero as $X$
tends to infinity.  However, this is not the case, and we will
show in fact that the conditions (1.1) in this statement are
optimal (see also [\Ro] for an announcement).

\proclaim Theorem 1.1. There exists a real number $\xi$ which is
neither rational nor quadratic irrational and which has the
property that, for a suitable constant $c>0$, the inequalities
(1.1) have a non-zero solution $(x_0,x_1,x_2)\in \bZ^3$ for any
real number $X\ge 1$.  Any such number is transcendental over
$\bQ$ and the set of these real numbers is countable.
\par

Let us call {\it extremal} any real number that has the property
stated in Theorem 1.1.  For such a number $\xi$, solutions to
(1.1) provide simultaneous approximations of $\xi$ and $\xi^2$ by
rational numbers $x_1/x_0$ and $x_2/x_0$ with the same denominator
(assuming $x_0\neq 0$).  The fact that, for any sufficiently large
$X$, these are sharper than expected by an application of
Dirichlet box principle (see for example Theorem 1A in Chapter II
of [\Sc]) reminds of Cassels' counter-example in questions of
algebraic independence (see Theorem XIV in Chapter V of [\Ca]).
Another situation, even closer to that of Cassels, where these
numbers satisfy estimates that are better than expected from the
box principle for all sufficiently large value of the parameter
$X$ is the following (compare with Theorem 1B in Chapter II of
[\Sc]).

\proclaim Theorem 1.2. Let $\xi$ be an extremal real number. There
exists a constant $c>0$ depending only on $\xi$ such that for any
real number $X\ge 1$ the inequalities
$$
|x_0| \le X,
  \quad
|x_1| \le X,
  \quad
|x_0\xi^2+x_1\xi+x_2| \le cX^{-\gamma^2} \eqno(1.3)
$$
have a non-zero solution $(x_0,x_1,x_2)\in \bZ^3$.
\par

Indeed, since $\gamma^2\simeq 2.618 >2$, the volume of the convex
body of $\bR^3$ defined by (1.3) tends to zero as $X$ tends to
infinity, a situation where we do not expect in principle to find
a non-zero solution for each sufficiently large $X$.

Examples of extremal real numbers are the numbers $\xi_{a,b}$
whose continued fraction expansion $\xi_{a,b} =
[0,a,b,a,a,b,a,b,\dots]$ is given by the Fibonacci word on two
distinct positive integers $a$ and $b$ (the infinite word
$abaabab\dots$ starting with $a$ which is the fixed point of the
substitution sending $a$ to $ab$ and $b$ to $a$).  We don't know
if, like these numbers $\xi_{a,b}$, all extremal real numbers have
bounded partial quotients but they satisfy the following weaker
measure of approximation by rational numbers.

\proclaim Theorem 1.3. Let $\xi$ be an extremal real number. Then,
there exist constants $c>0$ and $t>0$ such that, for any rational
number $\alpha\in\bQ$, we have
$$
| \xi-\alpha | \ge c H(\alpha)^{-2} (1+\log H(\alpha))^{-t}.
$$
\par

The above mentioned Fibonacci continued fractions $\xi_{a,b}$ are
an example of the more general Sturmian continued fractions
studied by J.-P.~Allouche, J.~L.~Davison, M.~Queff\'elec and
L.~Q.~Zamboni in [\ADQZ].  These authors showed that Sturmian
continued fractions are transcendental by showing that they admit
very good approximations by quadratic real numbers.  In the
present situation, we get an optimal measure of approximation by
such numbers.

\proclaim Theorem 1.4. Let $\xi$ be an extremal real number. There
is a value of $c>0$ such that the inequality
     $$ |\xi-\alpha| \le c H(\alpha)^{-2\gamma^2} $$
has infinitely many solutions in algebraic numbers $\alpha\in\bC$
of degree $2$ and another value of $c>0$ such that the same
inequality has no solution in algebraic numbers $\alpha\in\bC$ of
degree at most $2$. \par

It is generally believed that, for any real number $\xi$ that is
not algebraic over $\bQ$ of degree at most $3$ and for any
$\epsilon>0$, there should exist infinitely many algebraic
integers $\alpha$ of degree at most $3$ which satisfy
$|\xi-\alpha| \le H(\alpha)^{-3+\epsilon}$ (see page 259 of
[\Sc]). A recent result of Y.~Bugeaud and O.~Teuli\'e shows that
this is true of almost all real numbers $\xi$ in the sense of
Lebesgue measure (see Corollary 1 to Theorem 6 of [\BT]).  By
Schmidt's subspace theorem, this is also true of all algebraic
numbers of degree $>3$. However, we will present a criterion which
suggests that this is probably not true for extremal real numbers
and that, for these numbers, the correct exponent of approximation
is probably $\gamma^2$ instead of $3$ \footnote{($*$)}{Since this
paper was written the author found a class of extremal real
numbers for which this exponent of approximation is indeed
$\gamma^2$; see {\it Approximation to real numbers by cubic
algebraic integers II}, Annals of Math.\ (to appear),
arXiv:math.NT/0210182.}. As an application of this criterion, we
will derive the following measure of approximation.

\proclaim Theorem 1.5. Let $\xi$ be an extremal real number. There
exists a constant $c>0$ such that, for any algebraic integer
$\alpha\in\bC$ of degree $\le 3$, we have
    $$ |\xi-\alpha|\ge c H(\alpha)^{-\gamma^2-1}. $$
\par

This does not disprove the natural conjecture stated above, but it
sheds doubts on it. By extension, it also sheds doubts on the
conjectures of Wirsing and Schmidt concerning approximation to a
real number by algebraic numbers of bounded degree, although
Davenport and Schmidt have established these conjectures in the
case of approximation by quadratic irrationals in [\DSa].

This paper is organized as follows.  Section 5 presents a
criterion for a real number to be extremal.  Its proof is based on
preliminary considerations from Sections 2, 3 and 4.  Section 6
provides examples of extremal real numbers which include the
Fibonacci continued fractions.  Sections 7, 8 and 9 are devoted to
measures of approximations of extremal real numbers by rational
numbers, by quadratic real numbers and by algebraic integers of
degree at most 3 respectively.  Section 9 also provides a
criterion for approximation by algebraic integers of degree $\le
3$, and Section 10 proves a partial converse to this criterion.

As suggested by the referee, the reader may, at a first reading,
skip Section 2 and part (iii) of Lemma 3.1 and go directly from
Theorem 5.1 to Section 7. He may also want to devise an
independent proof of Proposition 6.3 showing that Fibonacci
continued fractions are extremal real numbers, or look at the
short note [\Ro] for such a proof.  However, in a more linear
reading, the relevance of the algebraic tool developed in Section
2 will appear clearly in Corollary 5.2, a result which the author
sees as a key for further study of extremal real numbers and a
motivation for the constructions of Section 6.

\medskip
\noindent{\bf Notation.}  Given a fixed real number $\xi$, we
often write estimates of the form $A\ll B$ or $B\gg A$, where $A$
and $B$ are non-negative variable real numbers, to mean $A\le cB$
for some constant $c>0$ depending only on $\xi$. We also write
$A\sim B$ to mean $A\ll B\ll A$.

%
%

\section{2. A particular construction}

Let $A$ denote a commutative ring with unit (in practice, we will
take $A=\bZ$). Identifying any point $\ux=(x_0,x_1,x_2)\in A^3$
with the symmetric matrix
$$
\ux=\pmatrix{x_0 &x_1\cr x_1 &x_2\cr},
$$
we can talk of its determinant
$$
\det(\ux)=x_0x_2-x_1^2.
$$
Put
$$
J=\pmatrix{0 &1\cr -1 &0\cr}.
$$
Then, given any three points $\ux, \uy, \uz \in A^3$, we find,
with the natural notation for their coordinates,
$$
\eqalign{ -\ux J \uz J \uy
 &=
\pmatrix{x_0 &x_1\cr x_1 &x_2\cr} \pmatrix{z_2 &-z_1\cr -z_1
&z_0\cr} \pmatrix{y_0 &y_1\cr y_1 &y_2\cr}  \cr \cr
 &=
\pmatrix{x_0 &x_1\cr x_1 &x_2\cr} \pmatrix{
  \deter{y_0 &y_1\cr z_1 &z_2\cr}
  &\deter{y_1 &y_2\cr z_1 &z_2\cr} \cr \cr \cr
  -\deter{y_0 &y_1\cr z_0 &z_1\cr}
  &-\deter{y_1 &y_2\cr z_0 &z_1\cr} \cr}  \cr \cr
 &=
\pmatrix{
 \deter{x_0 &x_1\cr\cr \deter{y_0 &y_1\cr z_0 &z_1\cr}
                     &\deter{y_0 &y_1\cr z_1 &z_2\cr}\cr}
 &\deter{x_0 &x_1\cr\cr \deter{y_1 &y_2\cr z_0 &z_1\cr}
                     &\deter{y_1 &y_2\cr z_1 &z_2\cr}\cr} \cr \cr \cr
 \deter{x_1 &x_2\cr\cr \deter{y_0 &y_1\cr z_0 &z_1\cr}
                     &\deter{y_0 &y_1\cr z_1 &z_2\cr}\cr}
 &\deter{x_1 &x_2\cr\cr \deter{y_1 &y_2\cr z_0 &z_1\cr}
                     &\deter{y_1 &y_2\cr z_1 &z_2\cr}\cr} \cr}.\cr}
$$
Moreover, if, in this matrix, we subtract from the element of row
1 and column 2 the element of row 2 and column 1, we find a
polynomial in $\ux$, $\uy$ and $\uz$ which happens to be
$$
\trace(J\ux J\uz J\uy) = \det(\ux,\uy,\uz).  \eqno(2.1)
$$
So, the product $-\ux J\uz J\uy$ is symmetric if and only if
$\det(\ux,\uy,\uz)=0$.    In this case, we denote by
$[\ux,\uy,\uz]$ the corresponding element of $A^3$:
$$
[\ux,\uy,\uz]=-\ux J\uz J\uy. \eqno(2.2)
$$

\proclaim Lemma 2.1. For any $\ux,\uy,\uz\in A^3$ with
$\det(\ux,\uy,\uz)=0$ and any $\uw\in A^3$, we have:
\item{(i)} $\det [\ux,\uy,\uz] = \det(\ux)\det(\uy)\det(\uz)$,
\item{(ii)} $\det(\uw,\uy,[\ux,\uy,\uz]) = \det(\uy)\det(\uw,\uz,\ux)$,
\item{(iii)} $\det(\ux,\uy,[\ux,\uy,\uz])=0$ and
             $[\ux,\uy,[\ux,\uy,\uz]] = \det(\ux)\det(\uy)\uz$.
\par

\proof The formula (i) follows from (2.2) since $\det(J)=1$.  To
prove (ii) and (iii), we note that, for any $\uw\in A^3$, we have
$$
\uw J\uw J = J\uw J\uw = -\det(\uw)I
$$
where $I$ denotes the identity matrix. Combining this observation
with (2.1), we find
$$
\eqalign{
\det(\uw,\uy,[\ux,\uy,\uz])
  &=\trace(J\uw J(-\ux J\uz J\uy)J\uy) \cr
  &=\det(\uy) \trace(J\uw J\ux J\uz) \cr
  &=\det(\uy) \det(\uw,\uz,\ux). \cr}
$$
This proves (ii) and the first half of (iii).  For the second
half, we compute
$$
-\ux J(-\ux J \uz J\uy) J \uy
  = (\ux J\ux J) \uz (J\uy J\uy)
  = \det(\ux)\det(\uy)\uz.
$$

%
%

\section{3. Some estimates}

Let $\xi$ be any real number.  For each integral point
$\ux=(x_0,x_1,x_2)\in \bZ^3$, we define
$$
\|\ux\|=\max\{|x_0|,|x_1|,|x_2|\}
  \and
L(\ux) = \max\{|x_0\xi-x_1|, |x_0\xi^2-x_2|\}.
$$
For a square matrix $A$, we also denote by $\|A\|$ the absolute
value of its determinant.  In the following lemma, we recall
estimates for determinants that appear in the study of Davenport
and Schmidt in Section 3 of [\DSb], and we provide estimates
relative to the construction of the previous paragraph.

\proclaim Lemma 3.1. Let $\ux,\uy,\uz\in\bZ^3$.
\item{(i)} For any choice of integers $r,s,t,u \in \{0,1,2\}$
with $s-r=u-t$, we have
$$
\left\|\matrix{x_r &x_s\cr y_t &y_u\cr}\right\|
  \ll \|\ux\| L(\uy) + \|\uy\| L(\ux).
$$
\item{(ii)} We have
$ |\det(\ux,\uy,\uz)|
  \ll \|\ux\| L(\uy) L(\uz) + \|\uy\| L(\ux) L(\uz)
      + \|\uz\| L(\ux)L(\uy).
$
\item{(iii)} Put $\uw=[\ux,\ux,\uy]$. Then, we have
$$
\|\uw\| \ll \|\ux\|^2 L(\uy) + \|\uy\| L(\ux)^2
  \and
L(\uw) \ll \big(\|\ux\| L(\uy) + \|\uy\| L(\ux)\big)L(\ux).
$$
\par

\proof  All these estimates are based on the multilinearity of the
determinant.  We simply prove (iii) since (i) and (ii) follow from
the computations in the proofs of Lemmas 3 and 4 of [\DSb]. Write
$\uw=(w_0,w_1,w_2)$.  Since
$$
w_0 = \deter{x_0 &x_1 \cr\cr \deter{x_0 &x_1\cr y_0 &y_1\cr}
            &\deter{x_0 &x_1\cr y_1 &y_2\cr}\cr}
    = \deter{x_0 &x_1-x_0\xi \cr\cr \deter{x_0 &x_1\cr y_0 &y_1\cr}
            &\deter{x_0 &x_1\cr y_1-y_0\xi &y_2-y_1\xi\cr}\cr},
$$
we obtain, using part (i) and noting that $L(\ux)\ll \|\ux\|$,
$$
|w_0|
  \ll \|\ux\| \absdeter{x_0 &x_1\cr y_1-y_0\xi &y_2-y_1\xi\cr}
      + L(\ux) \absdeter{x_0 &x_1\cr y_0 &y_1\cr}
  \ll \|\ux\|^2 L(\uy) + \|\uy\| L(\ux)^2.
$$
Similarly, we find, for $j=0,1$,
$$
|w_{j+1}-w_j\xi|
  = \absdeter{x_1-x_0\xi &x_2-x_1\xi \cr\cr
           \deter{x_j &x_{j+1}\cr y_0 &y_1\cr}
          &\deter{x_j &x_{j+1}\cr y_1 &y_2\cr}\cr}
  \ll L(\ux) (\|\ux\| L(\uy) + \|\uy\| L(\ux))
$$
and the conclusion follows.

%
%

\section{4. The sequence of minimal points}

Let $\xi$ be a real number which is neither rational nor quadratic
over $\bQ$.  Then, the function $L\colon \bZ^3\to\bR$ defined in
Section 3 has the property that, for $\ux,\uy\in\bZ^3$ with
$x_0,y_0\neq 0$, we have $L(\ux)=L(\uy)$ if and only if
$\ux=\pm\uy$. Accordingly, for any real number $X\ge 1$, the
finite set of points $\ux=(x_0,x_1,x_2)$ of $\bZ^3$ satisfying
$$
1\le x_0\le X
  \and
L(\ux) = \max\{|x_0\xi-x_1|, |x_0\xi^2-x_2|\} \le 1
$$
contains exactly one point for which $L(\ux)$ is minimal.
Following Davenport and Schmidt, we call it the {\it minimal point
corresponding to $X$}.

Clearly, if $\ux$ is a minimal point for some real number $X\ge
1$, then it is a minimal point for $x_0$.  So, we may order the
minimal points according to their first coordinates.  Note also
that the coordinates of a minimal point are relatively prime as a
set. Therefore, distinct minimal points are linearly independent
over $\bQ$ as elements of $\bQ^3$.  In particular, any two
consecutive minimal points are linearly independent over $\bQ$.

To state the next lemma, we define the height of a $2\times 3$
matrix $A$, denoted $\|A\|$, to be the maximal absolute value of
its minors of order $2$, and we define the height of a
$2$-dimensional subspace $V$ of $\bQ^3$, denoted $H(V)$, to be the
height of any $2\times 3$ matrix whose rows form a basis of
$V\cap\bZ^3$.

\proclaim Lemma 4.1.  Let $\ux$ be a minimal point, let $\uy$ be
the next minimal point and let $V=\langle \ux,\uy \rangle_\bQ$ be
the subspace of $\bQ^3$ generated by $\ux$ and $\uy$.  Then,
$\{\ux,\uy\}$ is a basis of $V\cap\bZ^3$ over $\bZ$ and we have
$H(V)\sim \|\uy\| L(\ux)$.
\par

\proof Write $\ux=(x_0,x_1,x_2)$ and $\uy=(y_0,y_1,y_2)$.  Arguing
like Davenport and Schmidt in their proof of Lemma 2 of [\DSa], we
observe that, if $\{\ux,\uy\}$ were not a basis of $V\cap\bZ^3$,
there would exist a point $\uz=(z_0,z_1,z_2)$ of $\bZ^3$ with
$z_0>0$, which could be written in the form $\uz=r\ux+s\uy$ for
rational numbers $r$ and $s$ with $\max\{|r|,|s|\} \le 1/2$. We
would then get
$$
z_0 \le |r|x_0 + |s|y_0 < y_0
  \and
L(\uz) \le |r| L(\ux) + |s| L(\uy) < L(\ux),
$$
in contradiction with the fact that $\uy$ is the next minimal
point after $\ux$.  Thus, $\{\ux,\uy\}$ is a basis of $V$ and,
using Lemma 3.1 (i), we find
$$
H(V) = \absdeter{x_0 &x_1 &x_2\cr y_0 &y_1 &y_2\cr}
  \ll \|\uy\| L(\ux) + \|\ux\| L(\uy) \ll \|\uy\| L(\ux).
$$
To derive lower bounds for $H(V)$, put $u_i=x_i-x_0\xi^i$ and
$v_i=y_i-y_0\xi^i$ for $i=1,2$, and choose an index $j$ for which
$|u_j|=L(\ux)$.  We find
$$
H(V) \ge \absdeter{x_0 &x_j\cr y_0 &y_j\cr}
  = |x_0 v_j - y_0 u_j|
  \ge y_0 |u_j| - x_0 |v_j|
  \ge (y_0-x_0) L(\ux).
$$
Since the point $\uz=\uy-\ux$ has its first coordinate
$z_0=y_0-x_0$ satisfying $1\le z_0 <y_0$, and since $\uy$ is the
next minimal point after $\ux$, there must exist an index $i$ for
which $|z_i-z_0\xi^i| = |v_i-u_i| > L(\ux)$. Thus, we also find
$$
H(V)
  \ge |x_0 v_i - y_0 u_i|
  \ge |x_0(v_i-u_i) - (y_0-x_0)u_i|
  \ge x_0 L(\ux) - (y_0-x_0) L(\ux).
$$
Combining the two preceding displayed inequalities, we find
$$
3H(V) \ge x_0 L(\ux) + (y_0-x_0) L(\ux) \ge y_0 L(\ux),
$$
and therefore $H(V) \gg \|\uy\| L(\ux)$.

%
%

\section{5. The set of extremal numbers}

Recall that, in the introduction, we defined a real number $\xi$
to be {\it extremal\/} if it is not rational nor quadratic
irrational and if there exists a constant $c>0$ such that the
inequalities
$$
1\le |x_0|\le X,
  \quad
|x_0\xi-x_1| \le c X^{-1/\gamma},
  \quad
|x_0\xi^2-x_2| \le c X^{-1/\gamma},
  \eqno(5.1)
$$
have a solution in integers $x_0,x_1,x_2$ for any real number
$X\ge 1$.  By virtue of the subspace theorem of Schmidt, such a
real number is transcendental over $\bQ$ (see for example Theorem
1B in Chapter VI of [\Sc]).

In this section, we characterize the extremal real numbers by a
stronger approximation property and show that the corresponding
approximation triples satisfy a certain recurrence relation
involving the operation $[\ux,\uy,\uz]$ defined in Section 2.  We
conclude that the set of extremal numbers is at most countable.

\proclaim Theorem 5.1.  A real number $\xi$ is extremal if and
only if there exists an increasing sequence of positive integers
$(Y_k)_{k\ge 1}$ and a sequence of points $(\uy_k)_{k\ge 1}$ of
$\bZ^3$ such that, for all $k\ge 1$, we have
$$
Y_{k+1}  \sim Y_k^\gamma,
  \quad
\|\uy_k\| \sim Y_k,
  \quad
L(\uy_k) \sim Y_k^{-1},
$$
$$
1\le |\det(\uy_k)| \ll 1
  \and
1\le |\det(\uy_k,\uy_{k+1},\uy_{k+2})| \ll 1,
$$
where the norm $\|\ \|$ and the map $L=L_\xi$ are as defined in
Section 3.
\par

\proof Assume first that there exist such sequences $(\uy_k)_{k\ge
1}$ and $(Y_k)_{k\ge 1}$ for a given real number $\xi$.  Then,
$\xi$ is not rational nor quadratic irrational because otherwise
there would exist integers $p,q,r$ not all zero such that $p + q
\xi + r \xi^2 = 0$ and consequently, for any point
$\uy_k=(y_{k,0}, y_{k,1}, y_{k,2})$, we would have
$$
|p y_{k,0} + q y_{k,1} + r y_{k,2}|
  = |q(y_{k,1}-y_{k,0}\xi) + r(y_{k,2}-y_{k,0}\xi^2)|
  \ll Y_k^{-1}.
$$
Since the left hand side of this inequality is an integer, this
integer would therefore be zero for all sufficiently large $k$, in
contradiction with the fact that the determinant of any three
consecutive points of the sequence $(\uy_k)_{k\ge 1}$ is non-zero.
Moreover, for any sufficiently large real number $X$, there exists
an index $k\ge 1$ such that $Y_k \le X < Y_{k+1}$ and then the
point $\uy_k$ satisfies
$$
\|\uy_k\| \ll Y_k \ll X
  \and
L(\uy_k) \ll Y_k^{-1} \ll Y_{k+1}^{-1/\gamma}
         \ll X^{-1/\gamma}.
$$
Therefore, $\xi$ is extremal.

Conversely, assume that $\xi$ is an extremal real number and
consider the sequence $(\ux_i)_{i\ge 1}$ of minimal points of
$\bZ^3$ associated to $\xi$, ordered according to their first
coordinate (see Section 4 above).  For each $i\ge 1$, denote by
$X_i$ the first coordinate of $\ux_i$ and put $L_i=L(\ux_i)$.
Choose also $c>0$ such that (5.1) has an integral solution for
each $X\ge 1$.  Then, by the definition of minimal points, we have
$$
L_i \le c X_{i+1}^{-1/\gamma},
  \quad (i\ge 1)
  \eqno(5.2)
$$
(see formula (30) in [\DSb]).  By Lemma 2 of [\DSb], there is an
index $i_0\ge 2$ such that $\det(\ux_i)\neq 0$ for all $i\ge i_0$.
Moreover, Lemma 5 of [\DSb] shows that there are infinitely many
indices $i> i_0$ such that $\ux_{i-1}$, $\ux_i$ and $\ux_{i+1}$
are linearly independent over $\bQ$.  Let $I$ denote the set of
all such integers $i$ and, for each integer $k\ge 1$, denote by
$i_k$ the $k$-th element of $I$.  We claim that the sequences
$(\uy_k)_{k\ge 1}$ and $(Y_k)_{k\ge 1}$ given by
$$
\uy_k = \ux_{i_k} \and Y_k=X_{i_k}, \quad (k\ge 1),
$$
have the required properties.

To prove this, fix any index $k\ge 1$ and put $i=i_k$.  The
inequality (5.2) gives
$$
L_{i-1} \ll X_i^{-1/\gamma} \and L_i \ll X_{i+1}^{-1/\gamma}.
$$
Then, by part (i) of Lemma 3.1 and the fact that $\det(\ux_i)$ is
a non-zero integer, we find
$$
1\le |\det(\bx_i)| \ll X_iL_i \ll X_i X_{i+1}^{-1/\gamma}
\eqno(5.3)
$$
and so $X_{i+1} \ll X_i^\gamma$.  Similarly, using part (ii) of
Lemma 3.1 and the fact that the integral points $\ux_{i-1}$,
$\ux_i$ and $\ux_{i+1}$ are linearly independent over $\bQ$, we
then find
$$
1 \le |\det(\bx_{i-1}, \bx_i, \bx_{i+1})|
  \ll X_{i+1} L_i L_{i-1}
  \ll X_{i+1}^{1/\gamma^2} L_{i-1}
  \ll X_{i+1}^{1/\gamma^2} X_i^{-1/\gamma}
  \ll 1.
$$
From these last estimates, we deduce
$$
X_{i+1} \sim X_i^\gamma,
  \quad
L_{i-1} \sim X_i^{-1/\gamma}
  \and
L_i \sim X_{i+1}^{-1/\gamma} \sim X_i^{-1}. \eqno(5.4)
$$
In particular, since $\ux_i=\uy_k$ and $X_i=Y_k$, we obtain, from
(5.3) and (5.4),
$$
\|\uy_k\| \sim Y_k,
  \quad
L(\uy_k) \sim Y_k^{-1}
  \and
1\le |\det(\uy_k)| \ll 1.
  \eqno(5.5)
$$

Now, let $V=\langle \ux_i,\ux_{i+1} \rangle_\bQ$ and let $j\ge
i+1$ be the largest integer such that $V$ contains $\ux_i,
\ux_{i+1}, \dots, \ux_j$. Then we have
$$
V=\langle \ux_i,\ux_{i+1} \rangle_\bQ
 =\langle \ux_{j-1},\ux_j \rangle_\bQ
\eqno(5.6)
$$
and, since $\ux_{j+1} \notin V$, we deduce that the points
$\ux_{j-1}$, $\ux_j$ and $\ux_{j+1}$ are linearly independent over
$\bQ$.  Moreover, since $\ux_i,\dots,\ux_j\in V$, the index $j$ is
the smallest integer $>i$ with the latter property.  Therefore, we
have $j=i_{k+1}$.  Consequently the estimates (5.4) also apply
with $i$ replaced by $j$. Lemma 4.1 then gives
$$
H\big( \langle \ux_i,\ux_{i+1} \rangle_\bQ \big)
  \sim X_{i+1} L_i
  \sim X_i^{1/\gamma}
\and H\big( \langle \ux_{j-1},\ux_j \rangle_\bQ \big)
  \sim X_j L_{j-1}
  \sim X_j^{1/\gamma^2}.
$$
By virtue of (5.6), this implies $X_j \sim X_i^\gamma$ and, since
$X_i=Y_k$ and $X_j=Y_{k+1}$, we conclude that
$$
Y_{k+1}\sim Y_k^\gamma.
  \eqno(5.7)
$$

The preceding argument also shows that $\langle \uy_k,\uy_{k+1}
\rangle_\bQ$ contains $\ux_{i+1}$ and $\ux_{j-1}$. Therefore,
$\langle \uy_{k+1},\uy_{k+2} \rangle_\bQ$ contains $\ux_{j+1}$ and
so $\langle \uy_k,\uy_{k+1},\uy_{k+2} \rangle_\bQ$ contains the
three linearly independent points $\ux_{j-1}$, $\ux_j$ and
$\ux_{j+1}$.  This implies that $\uy_k$, $\uy_{k+1}$ and
$\uy_{k+2}$ are linearly independent over $\bQ$ and, combining
Lemma 3.1 (ii) with (5.5) and (5.7), we find
$$
1\le |\det(\uy_k,\uy_{k+1},\uy_{k+2})|
 \ll \|\uy_{k+2}\| L(\uy_k) L(\uy_{k+1})
 \ll Y_{k+2} Y_k^{-1} Y_{k+1}^{-1}
 \ll 1.
$$
This estimate together with (5.5) and (5.7) proves that the
sequences $(\uy_k)_{k\ge 1}$ and $(Y_k)_{k\ge 1}$ have all the
required properties. \cqfd

\proclaim Corollary 5.2. Let $\xi$ be an extremal real number and
let $(\uy_k)_{k\ge 1}$ and $(Y_k)_{k\ge 1}$ be as in the statement
of Theorem 5.1.  Then, for any sufficiently large integer $k\ge
3$, the point $\uy_{k+1}$ is a non-zero rational multiple of
$[\uy_k,\uy_k,\uy_{k-2}]$.
\par

\proof Let $k\ge 4$ be an integer and let $\uw = [\uy_k, \uy_k,
\uy_{k+1}]$. By Lemma 2.1 (i), we have
$$
\det(\uw) = \det(\uy_k)^2 \det(\uy_{k+1}).
$$
Thus $\det(\uw)\neq 0$ and consequently, $\uw$ is a non-zero point
of $\bZ^3$.  Using Lemma 3.1 (iii), we find
$$
\|\uw\| \ll Y_k^2Y_{k+1}^{-1} \sim Y_{k-2}
  \and
L(\uw) \ll Y_{k+1} Y_{k}^{-2} \sim Y_{k-2}^{-1}.
$$
Then, applying Lemma 3.1 (ii), we obtain
$$
\eqalign{
|\det(\uw,\uy_{k-3},\uy_{k-2})|
  &\ll Y_{k-2} Y_{k-2}^{-1} Y_{k-3}^{-1} \sim Y_{k-3}^{-1}, \cr
|\det(\uw,\uy_{k-2},\uy_{k-1})|
  &\ll Y_{k-1} Y_{k-2}^{-2} \sim Y_{k-3}^{-1/\gamma}. \cr}
$$
But the above two determinants are integers. So, if $k$ is
sufficiently large, they must be zero.  Since $\uy_{k-3}$,
$\uy_{k-2}$ and $\uy_{k-1}$ are linearly independent over $\bQ$,
this then forces $\uw$ to be a rational multiple of $\uy_{k-2}$.
Since, by Lemma 2.1 (iii), we have
$$
[\uy_k,\uy_k,\uw] = \det(\uy_k)^2 \uy_{k+1},
$$
we conclude that $\uy_{k+1}$ is a rational multiple of
$[\uy_k,\uy_k,\uy_{k-2}]$. \cqfd

\medskip
As we will see in the next section, this corollary provides a way
to construct extremal numbers.  It also provides the following
formula.

\proclaim Corollary 5.3. In the notation of Corollary 5.2, we have
$$
\det(\uy_{k-2},\uy_{k-1},\uy_k)\,\uy_{k+1}
 = \det(\uy_{k-2},\uy_{k-1},\uy_{k+1})\,\uy_{k}
   + \det(\uy_{k-1},\uy_{k},\uy_{k+1})\,\uy_{k-2}
$$
for any sufficiently large value of $k$.
\par

\proof By Lemma 2.1 (ii), we have $\det(\uy_{k-2}, \uy_k, [\uy_k,
\uy_k, \uy_{k-2}])=0$ and this, by Corollary 5.2, implies that the
points $\uy_{k-2}$, $\uy_k$ and $\uy_{k+1}$ are linearly dependent
over $\bQ$.  Writing $\uy_{k+1}=a\uy_k+b\uy_{k-2}$ with unknown
rational coefficients $a$ and $b$ and using the multilinearity of
the determinant, we get $\det(\uy_{k-2}, \uy_{k-1}, \uy_{k+1})= a
\det(\uy_{k-2}, \uy_{k-1}, \uy_{k})$ and $\det(\uy_{k-1}, \uy_{k},
\uy_{k+1})= b \det(\uy_{k-1}, \uy_{k}, \uy_{k-2})$ and the formula
follows. \cqfd

\medskip
\proclaim Corollary 5.4. The set $\cE$ of extremal real numbers is
at most countable.
\par

\proof We prove this by constructing an injective map $\varphi
\colon\cE\to(\bZ^3)^3$ as follows. For each extremal real number
$\xi\in\cE$, we choose a corresponding sequence $(\uy_k)_{k\ge 1}$
as in Theorem 5.1. Then, in accordance with Corollary 5.3, we
choose an index $i\ge 3$ such that, for each $k\ge i$, the point
$\uy_{k+1}$ is a non-zero rational multiple of $[\uy_k, \uy_k,
\uy_{k-2}]$, and we put $\varphi(\xi)=(\uy_{i-2}, \uy_{i-1},
\uy_i)$.  Since the knowledge of $\uy_{i-2}$, $\uy_{i-1}$, $\uy_i$
determines all points $\uy_k$ with $k\ge i-2$ up to a non-zero
rational multiple, it determines uniquely the corresponding ratios
$y_{k,1}/y_{k,0}$ and also their limit $\xi = \lim_{k\to\infty}
y_{k,1}/y_{k,0}$.  Thus, $\varphi$ is injective. \cqfd

\medskip
In order to complete the proof of Theorem 1.1, it remains to show
that the set of extremal numbers is infinite.  This is the object
of the next section.

%
%

\section{6. Construction of extremal numbers}

As in Section 2, we identify each point $(y_0,y_1,y_2)$ of $\bZ^3$
with the corresponding symmetric matrix with integral coefficients
$\pmatrix{y_0 &y_1\cr y_1 &y_2\cr}$.  Before stating the main
result of this section, we first prove a lemma.

\proclaim Lemma 6.1. Let $A$, $B$ be non-commuting symmetric
matrices in $\GL_2(\bZ)$.  Consider the sequence $(\uy_k)_{k\ge
-1}$ constructed recursively by putting
$$
\uy_{-1}=B^{-1},
  \quad
\uy_0=I,
  \quad
\uy_1=A
  \and
\uy_k = [\uy_{k-1}, \uy_{k-1}, \uy_{k-3}]
  \quad \hbox{for $k\ge 2$,}
  \eqno(6.1)
$$
where $I$ denotes the identity $2\times 2$ matrix.  Then, for each
integer $k\ge -1$, we have
$$
|\det(\uy_k)|= 1
  \and
|\det(\uy_k,\uy_{k+1},\uy_{k+2})| = |\trace(JAB)| \neq 0
  \eqno(6.2)
$$
where $J$ is as defined in Section 2.  Moreover, for $k\ge 1$, we
have the alternative recurrence relation
$$
\uy_k = \pm \uy_{k-1} S \uy_{k-2}
  \quad\hbox{where}\quad
S = \cases{AB & if $k$ is odd,\cr BA &if $k$ is even.\cr}
  \eqno(6.3)
$$
\par

\proof For $k\ge 2$, the recurrence relation in (6.1) implies,
using Lemma 2.1 (i) and (ii),
$$
\eqalign{ \det(\uy_k)
  &= \det(\uy_{k-1})^2\det(\uy_{k-3}),  \cr
\det(\uy_{k-2},\uy_{k-1},\uy_k)
  &= - \det(\uy_{k-1}) \det(\uy_{k-3}, \uy_{k-2}, \uy_{k-1}). \cr}
$$
Since $A$ and $B$ belong to $\GL_2(\bZ)$, we also have
$\det(\uy_k) = \pm 1$ for $k=-1,0,1$, and therefore we deduce, by
recurrence, that, for all $k\ge -1$, we have
$$
|\det(\uy_k)| = 1
  \and
|\det(\uy_k,\uy_{k+1},\uy_{k+2})| = |\det(\uy_{-1},\uy_0,\uy_1)|.
$$
Moreover, using (2.1), we find $\det(\uy_{-1},\uy_0,\uy_1) =
\det(B) \trace(JAB)\neq 0$ since $AB\neq BA$. This proves (6.2).

Similarly, a short computation shows that (6.3) holds for $k=1$
and, assuming that it holds for some integer $k\ge 1$, we find,
using (2.2),
$$
\uy_{k+1}
  = - \uy_k J \uy_{k-2} J \uy_k
  = \pm \uy_{k-1} S \uy_{k-2} J \uy_{k-2} J \uy_k
  = \pm \uy_{k-1} S \uy_k
$$
which, by taking transposes, implies $\uy_{k+1} = \pm \uy_k {^tS}
\uy_{k-1}$.  Thus (6.3) holds for all $k\ge 1$. \cqfd

\medskip
\proclaim Theorem 6.2. Let $A$, $B$ and the sequence
$(\uy_k)_{k\ge -1}$ be as in Lemma 6.1.  Assume that all
coefficients of $A$ are non-negative and that all those of $AB$
are positive. Then, upon writing $\uy_k = (y_{k,0}, y_{k,1},
y_{k,2})$, we have that $y_{k,0}$ is a non-zero integer for all
$k\ge 2$ and that the limit
$$
\xi =\lim_{k\to\infty} y_{k,1}/y_{k,0}
$$
exists and is an extremal real number.  Moreover, letting
$Y_k=\|\uy_k\|$ for all $k$, the sequences $(\uy_k)_{k\ge 3}$ and
$(Y_k)_{k\ge 3}$ share, with respect to $\xi$, all the properties
stated in Theorem 5.1.
\par

\proof  By virtue of (6.2), the sequence $(\uy_k)_{k\ge -1}$
satisfies the last two conditions of Theorem 5.1.  Since $\uy_0$
and $\uy_1$ have non-negative coefficients, since their rows and
columns are all non-zero, and since $AB$ and its transpose $BA$
have positive coefficients, the recurrence relation (6.3) implies
that, for $k\ge 2$, the point $\uy_k$ has non-zero coefficients,
all of the same sign.  Thus, the integer $Y_k = \|\uy_k\|$ is
positive for $k\ge 2$ and there exists a constant $c_1 > 1$ such
that, for all $k\ge 4$, we have
$$
Y_{k-2} Y_{k-1} < Y_k \le c_1 Y_{k-2} Y_{k-1}.
  \eqno(6.4)
$$
In particular, the sequence $(Y_k)_{k\ge 3}$ is unbounded and
monotone increasing.  Similarly, we find that $\lim_{k\to\infty}
y_{k,0} = \infty$.

For each $k\ge 2$, define $I_k$ to be the interval of $\bR$ with
end points $y_{k,1}/y_{k,0}$ and $y_{k,2}/y_{k,1}$.  For $k\ge 3$,
the relation (6.3) shows that there are positive integers
$r,s,t,u$ depending on $k$ such that
$$
\pmatrix{y_{k,0} &y_{k,1}\cr y_{k,1} &y_{k,2}\cr}
  = \pm \pmatrix{y_{k-1,0} &y_{k-1,1}\cr y_{k-1,1} &y_{k-1,2}\cr}
        \pmatrix{r &s\cr t &u\cr}
$$
and consequently we find
$$
{y_{k,1}\over y_{k,0}}
  = {r y_{k-1,1} + t y_{k-1,2} \over r y_{k-1,0} + t y_{k-1,1}}
  \in I_{k-1}
\and {y_{k,2}\over y_{k,1}}
  = {s y_{k-1,1} + u y_{k-1,2} \over s y_{k-1,0} + u y_{k-1,1}}
  \in I_{k-1}.
$$
Thus, the intervals $I_k$ with $k\ge 2$ form a non-increasing
sequence $I_2 \supseteq I_3 \supseteq I_4 \supseteq \cdots$ and,
since the length of $I_k$ is $1/(y_{k,0}y_{k,1})$ tending to zero
as $k$ tends to infinity, their intersection is reduced to one
point $\xi$ with $\xi>0$.  The ratios $y_{k,1}/y_{k,0}$ and
$y_{k,2}/y_{k,1}$ therefore converge to that number $\xi$ and,
since $\xi\neq 0$, we get
$$
|y_{k,0}| \sim |y_{k,1}| \sim |y_{k,2}| \sim Y_k.
$$
In particular, the length of $I_k$ satisfies $|I_k| =
(y_{k,0}y_{k,1})^{-1} \sim Y_k^{-2}$ and so we find
$$
L(\uy_k)
  \sim \max_{j=0,1} |y_{k,j}\xi-y_{k,j+1}|
  \sim Y_k \max_{j=0,1} \left|\xi-{y_{k,j+1}\over y_{k,j}}\right|
  \sim Y_k |I_k|
  \sim Y_k^{-1}.
$$
To conclude from Theorem 5.1 that $\xi$ is extremal, it remains
only to show that $Y_k\sim Y_{k-1}^\gamma$.  To this end, we
observe that, for $k\ge 4$, the real number $q_k=Y_k
Y_{k-1}^{-\gamma}$ satisfies, by (6.4),
$$
q_{k-1}^{-1/\gamma} \le q_k \le c_1 q_{k-1}^{-1/\gamma}.
$$
By recurrence on $k$, this implies that, for all $k\ge 3$, we have
$c_2^{-1/\gamma} \le q_k \le c_2$ for the constant $c_2 =
\max\{c_1^\gamma, q_3, q_3^{-\gamma}\}$ and so $c_2^{-1/\gamma}
Y_{k-1}^\gamma \le Y_k \le c_2 Y_{k-1}^\gamma$ as required. \cqfd

\bigskip
Let $E=\{a, b\}$ be an alphabet of two letters and let $E^*$
denote the monoid of words on $E$ for the product given by
concatenation. The {\it Fibonacci sequence} in $E^*$ is the
sequence of words $(w_i)_{i\ge 0}$ defined recursively by
$$
w_0=b, \quad w_1=a \and w_i=w_{i-1}w_{i-2}, \quad (i\ge 2)
$$
(see Example 1.3.6 of [\Lo]).  Since, for every $i\ge 1$, the word
$w_i$ is a prefix of $w_{i+1}$, this sequence converges to an
infinite word $w=abaabab\dots$ called the {\it Fibonacci word} on
$\{a,b\}$.  We now turn to Fibonacci continued fractions.

\proclaim Corollary 6.3. Let $a$ and $b$ be distinct positive
integers and let
$$
\xi_{a,b} = [0,a,b,a,a,b,a, \dots] = 1/(a+1/(b+...))
$$
be the real number whose sequence of partial quotients starts with
$0$ followed by the elements of the Fibonacci word on $\{a,b\}$.
Then $\xi_{a,b}$ is an extremal real number.  It is the special
case of the construction of Theorem 6.2 corresponding to the
choice of
$$
A = \pmatrix{a &1\cr 1 &0\cr}
  \and
B = \pmatrix{b &1\cr 1 &0\cr}.
$$
\par

\medskip
Note that, if $\xi$ is an extremal real number, then $1/\xi$ is
also extremal.  So, the corollary implies that $1/\xi_{a,b} =
[a,b,a,a,b,\dots]$ as well is extremal. More generally, one finds
that, if $\xi$ is extremal, then the ratio $(r\xi+s)/(t\xi+u)$ is
also extremal for any choice of rational numbers $r,s,t,u$ with
$ru-st\neq 0$.

\proof  Recall that, if a real number $\xi$ has an infinite
continued fraction expansion $\xi=[a_0, a_1, a_2, \dots]$ and if,
for each $i\ge 0$, we write its convergent $[a_0,a_1,\dots,a_i]$
in the form $p_i/q_i$ with relatively prime integers $p_i$ and
$q_i>0$, then the sequences $(p_i)_{i\ge 0}$ and $(q_i)_{i\ge 0}$
satisfy recurrence relations which, in matrix form, can be written
for $i\ge 2$ as
$$
\pmatrix{q_i &q_{i-1}\cr p_i &p_{i-1}\cr}
  = \pmatrix{q_{i-1} &q_{i-2}\cr p_{i-1} &p_{i-2}\cr}
    \pmatrix{a_i &1\cr 1 &0\cr}
$$
(see for example Chapter I of [\Sc]).  When $a_0=0$, this gives
$$
\pmatrix{q_i &q_{i-1}\cr p_i &p_{i-1}\cr}
  = \pmatrix{a_1 &1\cr 1 &0\cr} \pmatrix{a_2 &1\cr 1 &0\cr}
    \cdots \pmatrix{a_i &1\cr 1 &0\cr}.
$$

Now, consider the morphism of monoids $\varphi\colon E^* \to
\GL_2(\bZ)$ determined by the conditions
$$
\varphi(a) = A  \and \varphi(b) = B.
$$
For each integer $k\ge 0$, let $m_k$ denote the word $w_{k+2}$
minus its last two letters and let $\ux_k=\varphi(m_k)$.  Then, we
have $m_0=1$ (the empty word), $m_1=a$ and, for $k\ge 2$,
$$
m_k = m_{k-1} s m_{k-2}
  \quad \hbox{where} \quad
s = \cases{ab &if $k$ is odd,\cr
           ba &if $k$ is even.\cr}
$$
In terms of the matrices $\ux_k$, this implies $\ux_0=I$,
$\ux_1=A$ and, for $k\ge 2$,
$$
\ux_k = \ux_{k-1} S \ux_{k-2}
  \quad \hbox{where} \quad
S = \cases{AB &if $k$ is odd,\cr
           BA &if $k$ is even.\cr}
$$
Therefore, if $(\uy_k)_{k\ge -1}$ denotes the sequence furnished
by Lemma 6.1 for the above choice of $A$ and $B$, we deduce from
(6.3) that $\uy_k = \pm \ux_k$ for any $k\ge 0$.  In particular,
$\ux_k$ is a symmetric matrix and, with the usual notation for
coordinates, we have
$$
{y_{k,1}\over y_{k,0}} = {x_{k,1}\over x_{k,0}}
$$
for all $k\ge 1$.  On one hand, Theorem 5.1 tells us that the
above ratio converges to an extremal real number $\xi$ as $k$
tends to infinity. On the other hand, since $m_k$ is a prefix of
the Fibonacci word $w$ and since $\ux_k = \varphi(m_k)$, the
observation made at the beginning of the proof shows that the same
ratio is the convergent of $\xi_{a,b}$ whose partial quotients are
$0$ followed by the elements of $m_k$. Therefore $\xi_{a,b}=\xi$
is extremal. \cqfd

%
%

\section{7. Approximation by rational numbers}

In this section, we prove Theorem 1.3, as a consequence of the
following result.

\proclaim Proposition 7.1. Let $\xi$ be an extremal real number.
There exists a positive constant $s$ such that, for any
sufficiently large real number $X$, the first minimum of the
convex body
$$
\cC(X)=\{\; (x_0,x_1)\in\bR^2 \;;\; |x_0| \le X \and |x_0\xi-x_1|
\le X^{-1} \;\}
$$
is bounded below by $(\log X)^{-s}$. \par

\proof  Let $(Y_k)_{k\ge 1}$ and $(\uy_k)_{k\ge 1}$ be the
sequences given by Theorem 5.1.  Let $X$ be a real number with
$X\ge Y_2$, let $\lambda=\lambda(X)$ be the first minimum of
$\cC(X)$, and let $(x_0,x_1)$ be a point of $\bZ^2$ which realizes
this minimum, so that we have $|x_0|\le \lambda X$ and
$|x_0\xi-x_1| \le \lambda X^{-1}$.  We choose an index $k\ge 2$
such that $Y_k\le X\le Y_{k+1}$ and consider the point $(z_0,z_1)$
of $\bZ^2$ given by
$$
(z_0,z_1)= (x_0,x_1)J\uy_{k+1} J\uy_{k-1}
$$
in the notation of Section 2. Writing $\uy_{k-1} = (y_0^*, y_1^*,
y_2^*)$ and $\uy_{k+1} = (y_0', y_1', y_2')$, and using the
estimates of Theorem 5.1 together with Lemma 3.1 (i), we find
$$
\eqalign{ |z_0|
   &= \absdeter{x_0 &x_1\cr\cr
        \deter{y_0^* &y_1^*\cr y_0' &y_1'\cr}
       &\deter{y_0^* &y_1^*\cr y_1' &y_2'\cr}\cr}
    = \absdeter{x_0 &x_0\xi-x_1\cr\cr
        \deter{y_0^* &y_1^*\cr y_0' &y_1'\cr}
       &\deter{y_0^* &y_1^*\cr y_0'\xi-y_1' &y_1'\xi-y_2'\cr}\cr}\cr\cr
   &\ll {|x_0|\over Y_k} + Y_k |x_0\xi-x_1|
    \ll {\lambda X\over Y_k} \cr}
$$
and
$$
\eqalign{
|z_0\xi-z_1|
  &= \absdeter{x_0 &x_1\cr\cr
       \deter{y_0^*\xi-y_1^* &y_1^*\xi-y_2^*\cr y_0' &y_1'\cr}
      &\deter{y_0^*\xi-y_1^* &y_1^*\xi-y_2^*\cr y_1' &y_2'\cr}\cr}
       \cr\cr
  &= \absdeter{x_0 &x_0\xi-x_1\cr\cr
       \deter{y_0^*\xi-y_1^* &y_1^*\xi-y_2^*\cr y_0' &y_1'\cr}
      &\deter{y_0^*\xi-y_1^* &y_1^*\xi-y_2^*\cr
              y_0'\xi-y_1' &y_1'\xi-y_2'\cr}\cr}\cr\cr
  &\ll {|x_0|\over Y_{k-1} Y_{k+1}} + Y_k |x_0\xi-x_1|
   \ll {\lambda Y_k\over X}. \cr}
$$
Thus, if we put $Z=2X/Y_k$, these relations imply that the first
minimum $\lambda(Z)$ of $\cC(Z)$ satisfies $\lambda(Z) \le c
\lambda(X)$ for some constant $c>0$ which is independent of $k$
and $X$.  We also note that $Z\ll X^{1/\gamma^2}$ and thus $Z\le
X^{1/2}$ if $X$ is sufficiently large, say $X\ge X_0\ge Y_2$.
Then, by choosing $s>0$ with $2^{s-1}\ge c$, we get, for $X\ge
X_0$,
$$
\lambda(X)(\log X)^s \ge 2 \lambda(Z) (\log Z)^s
$$
where $Z$ is a real number in the interval  $2\le Z\le X^{1/2}$.
Since $\lambda(X)(\log X)^s$ is bounded below by some positive
constant for $2\le X\le X_0^{1/2}$, this shows that
$\lambda(X)(\log X)^s$ tends to infinity with $X$. \cqfd

\medskip
\noindent {\it Proof of Theorem 1.3.}  Let $s$ be as in the
proposition and assume that a rational number $\alpha=p/q$ with
denominator $q\ge 3$ satisfies $|\xi-p/q| \le q^{-2}(\log
q)^{-2s-2}$. Then, by putting $X=q(\log q)^{s+1}$ and
$\lambda=(\log q)^{-s-1}$, we get $|q|\le \lambda X$ and
$|q\xi-p|\le \lambda/X$ so that the first minimum of $\cC(X)$ is
$\le \lambda$.  By the proposition, this implies $\lambda \ge
(\log X)^{-s}$ if $q$ is sufficiently large.  This in turn forces
an upper bound on $q$ and therefore on $H(p/q)$.  Thus Theorem 1.3
holds with $t=2s+2$. \cqfd

%
%

\section{8. Approximation by quadratic real numbers}

In this section, we prove Theorems 1.2 and 1.4 as consequences of
the following result.

\proclaim Proposition 8.1.  Let $\xi$ be an extremal real number
and let the sequences $(Y_k)_{k\ge 1}$ and $(\uy_k)_{k\ge 1}$ be
as in Theorem 5.1.  For each $k\ge 1$, we define
$$
Q_k(T) = \deter{1 &T &T^2\cr y_{k,0} &y_{k,1} &y_{k,2}\cr
                y_{k+1,0} &y_{k+1,1} &y_{k+1,2}\cr}.
$$
Then, $Q_k$ is a polynomial with integral coefficients which, for
all sufficiently large values of $k$, has degree $2$ and satisfies
$$
H(Q_k) \sim |Q'_k(\xi)| \sim Y_{k-1}
  \and
|Q_k(\xi)| \sim Y_{k+2}^{-1}.
$$
\par

\proof By Lemma 3.1 (i), we have
$$
H(Q_k) = \absdeter{y_{k,0} &y_{k,1} &y_{k,2}\cr
                   y_{k+1,0} &y_{k+1,1} &y_{k+1,2}\cr}
  \ll Y_{k+1} L(\uy_{k}) + Y_{k} L(\uy_{k+1})
  \ll Y_{k-1}.
$$
All the remaining estimates assume that $k$ is sufficiently large.
Since
$$
Q_k'' = 2 \deter{y_{k,0} &y_{k,1}-y_{k,0}\xi\cr
                    y_{k+1,0} &y_{k+1,1}-y_{k+1,0}\xi\cr}
  = 2 y_{k+1,0}(y_{k,0}\xi-y_{k,1}) + \cO(Y_{k-1}^{-1}),
$$
we get $Q_k''\neq 0$ for $k\gg 1$ using the lower bound on
$|y_{k,0}\xi-y_{k,1}|$ provided by Theorem 1.3. Similarly, we have
$$
Q'_k(\xi)
  = - \deter{y_{k,0}
            &y_{k,2} - 2 y_{k,1} \xi + y_{k,0}\xi^2\cr
             y_{k+1,0}
            &y_{k+1,2} - 2 y_{k+1,1} \xi + y_{k+1,0} \xi^2\cr}.
\eqno(8.1)
$$
Note that, for any $j\ge 1$,
$$
y_{j,0}(y_{j,2} - 2 y_{j,1} \xi + y_{j,0}\xi^2)
  = \det(\uy_j) +(y_{j,1}-y_{j,0}\xi)^2
  = \det(\uy_j) + \cO(Y_j^{-2}),
$$
and therefore
$$
| y_{j,2} - 2 y_{j,1} \xi + y_{j,0}\xi^2 |
  \sim Y_j^{-1}.
$$
Using this to estimate (8.1), we find
$$
|Q'_k(\xi)|
  = |y_{k+1,0}|\,|y_{k,2}-2y_{k,1}\xi+y_{k,0}\xi^2|
    + \cO(Y_{k-1}^{-1})
  \sim Y_{k-1}
$$
and, since $H(Q_k)\ll Y_{k-1}$, we get $H(Q_k)\sim Y_{k-1}$.
Finally, write
$$
y_{k+2,0}(1,\xi,\xi^2) = \uy_{k+2} + \uz
$$
where $\uz\in\bR^3$ has norm $\|\uz\| \ll Y_{k+2}^{-1}$. By
multilinearity of the determinant, this gives
$$
y_{k+2,0} Q_k(\xi)
  = \det(\uy_{k+2},\uy_k,\uy_{k+1}) + \det(\uz,\uy_k,\uy_{k+1}).
$$
Since $$|\det(\uz,\uy_k,\uy_{k+1})| \le 3\|\uz\|H(Q_k) \ll
Y_{k-1}/Y_{k+2},$$ we conclude that
$$
|Q_k(\xi)|
  = |y_{k+2,0}|^{-1}
    \big( |\det(\uy_k,\uy_{k+1},\uy_{k+2})|
         + \cO(Y_{k-1}/Y_{k+2}) \big)
  \sim Y_{k+2}^{-1}. \hbox{\cqfd}
$$

\bigskip
\noindent
{\it Proof of Theorem 1.2.} For any real number
$X$ with $X\ge H(Q_1)$, we choose an index $k\ge 1$ such that
$H(Q_k)\le X \le H(Q_{k+1})$.  Then upon writing
$Q_k(T)=x_0T^2+x_1T+x_2$, we find that $(x_0,x_1,x_2)$ is a
non-zero integral point with $|x_0|\le X$, $|x_1|\le X$ and
$$
|x_0\xi^2+x_1\xi+x_2|
  = |Q_k(\xi)|
  \sim Y_{k+2}^{-1}
  \sim H(Q_{k+1})^{-\gamma^2}
  \ll X^{-\gamma^2}.
\hbox{\cqfd}
$$

\medskip
As for Theorem 1.4, it follows immediately from the following more
precise result.

\proclaim Theorem 8.2. Let $\xi$ be an extremal real number. Then
there exist an integer $k_0\ge 1$ and positive constants $c_1$,
$c_2$ such that, for any $k\ge k_0$, the polynomial $Q_k$ of Lemma
8.1 is irreducible over $\bQ$ of degree $2$ and admits exactly one
root $\alpha_k$ with
$$
c_1 H(\alpha_k)^{-2\gamma^2}
  \le |\xi-\alpha_k|
  \le c_2 H(\alpha_k)^{-2\gamma^2}.
\eqno(8.2)
$$
There also exists a constant $c_3>0$ such that, for any algebraic
number $\alpha\in\bC$ of degree $\le 2$ over $\bQ$, distinct from
all $\alpha_k$ with $k\ge k_0$, we have
$$
|\xi-\alpha| \ge c_3 H(\alpha)^{-4}.
$$
\par

\proof By Proposition 8.1, the polynomial $Q_k$ has degree $2$ for
$k\gg1$. For these values of $k$, we may factor it as $Q_k(T) =
a_k (T-\alpha_k) (T-\beta_k)$ with the roots $\alpha_k$ and
$\beta_k$ ordered so that $|\xi-\alpha_k| \le |\xi-\beta_k|$. Then
the logarithmic derivative of $Q_k$ satisfies
$$
{1\over |\xi-\alpha_k|}-{1\over |\xi-\beta_k|}
  \le {|Q_k'(\xi)|\over |Q_k(\xi)|}
  = \left|{1\over \xi-\alpha_k} + {1\over \xi-\beta_k} \right|
  \le {2\over |\xi-\alpha_k|}.
\eqno(8.3)
$$
Using the estimates of the proposition, this gives
$$
|\xi-\alpha_k|
  \le 2{|Q_k(\xi)|\over |Q_k'(\xi)|}
  \sim H(Q_k)^{-2\gamma^2}
  \ll H(\alpha_k)^{-2\gamma^2}.
\eqno(8.4)
$$
Since $H(Q_k)$ tends to infinity with $k$, we deduce that
$\alpha_k$ tends to $\xi$ as $k$ tends to infinity and, by
comparison with Theorem 1.3, that $\alpha_k$ is irrational for all
sufficiently large values of $k$.  Thus, for $k\gg 1$, the
polynomial $Q_k$ is irreducible over $\bQ$. Its discriminant
$a_k^2(\alpha_k-\beta_k)^2$ is then a non-zero integer, and so we
find
$$
|\alpha_k-\beta_k| \ge |a_k|^{-1} \ge H(Q_k)^{-1}.
$$
In particular, we have $|\alpha_k-\beta_k|\ge 3|\xi-\alpha_k|$ for
any sufficiently large value of $k$, and (8.3) becomes
$$
{1\over 2|\xi-\alpha_k|}
  \le {|Q_k'(\xi)|\over |Q_k(\xi)|}
  \le {2\over |\xi-\alpha_k|}.
$$
On the other hand, the greatest common divisor of the coefficients
of $Q_k$ divides the non-zero integer $\det(\uy_{k-1}, \uy_k,
\uy_{k+1})$ and, since the absolute value of this integer is
bounded above by a constant that is independent of $k$, we also
have $H(\alpha_k) \sim H(Q_k)$ for $k\gg 1$.  Then, going back to
(8.4), we get refined estimates of the form (8.2) valid say for
all $k\ge k_0$, where $c_1$, $c_2$ are positive constants and
$k_0\ge 2$ is a fixed integer.  This proves the first part of the
Theorem.

Now, let $\alpha\in\bC$ be an algebraic number of degree $\le 2$
with $\alpha\neq\alpha_k$ for any index $k\ge k_0$. For these
values of $k$, Liouville's inequality gives
$$
|\alpha-\alpha_k|
  \ge c_3 H(\alpha)^{-2}H(\alpha_k)^{-2}
$$
where $c_3>0$ is an absolute constant (see for example Lemma 3 of
[\LR]).  Choose $k$ to be the smallest integer $\ge k_0$
satisfying
$$
H(\alpha)
  \le \left({c_3\over 2c_2}\right)^{1/2} H(\alpha_k)^\gamma.
$$
This choice of $k$ ensures that
$$
|\xi-\alpha_k|
  \le c_2 H(\alpha_k)^{-2\gamma^2}
  \le {c_3\over 2} H(\alpha)^{-2}H(\alpha_k)^{-2}
  \le {1\over 2} |\alpha-\alpha_k|,
$$
and therefore that
$$
|\xi-\alpha|
  \ge {1\over 2} |\alpha-\alpha_k|
  \gg H(\alpha)^{-2}H(\alpha_k)^{-2}.
$$
If $k>k_0$, we also have $H(\alpha)\gg H(\alpha_{k-1})^\gamma \gg
H(\alpha_k)$, while if $k=k_0$, we have $H(\alpha)\ge 1 \gg
H(\alpha_k)$.  So, the above inequality leads to $|\xi-\alpha| \gg
H(\alpha)^{-4}$.  \cqfd

%
%

\section{9.  Approximation by cubic algebraic integers}

For a real number $x$, let $\dZ{x} = \min\{|x-p|\,;\,p\in\bZ\}$
denote the distance from $x$ to a nearest integer.  In this
section, we fix an extremal real number $\xi$ and corresponding
sequences $(Y_k)_{k\ge 1}$ and $(\uy_k)_{k\ge 1}$ satisfying the
conditions of Theorem 5.1. In particular, we have
$\dZ{y_{k,0}\xi^j} \ll Y_k^{-1}$ for $k\ge 1$ and $j=0,1,2$. Here,
we prove a lower bound for $\dZ{y_{k,0}\xi^3}$ which implies the
measure of approximation of Theorem 1.5 through the following
proposition.

\proclaim Proposition 9.1.  Assume that there exist real numbers
$\delta$ and $c_1$ with $0\le \delta<1$ and $c_1>0$ such that
$$
\dZ{y_{k,0}\xi^3} \ge c_1 Y_k^{-\delta} \eqno(9.1)
$$
for all $k\ge 1$.  Then there exists a constant $c_2>0$ such that,
for any algebraic integer $\alpha\in\bC$ of degree $\le 3$, we
have
$$
|\xi-\alpha| \ge c_2 H(\alpha)^{-\theta}
$$
where $\theta=(\gamma^2+\delta/\gamma)/(1-\delta)$.
\par

\proof Let $\alpha$ be an algebraic integer of degree $\le 3$ and
let $P(T)=T^3+pT^2+qT+r$ be the product of its irreducible
polynomial over $\bZ$ by the appropriate power of $T$ which makes
it of degree $3$.  Since $H(P)=H(\alpha)$ we find, for all $k\ge
1$,
$$
\eqalign{ \dZ{y_{k,0}\xi^3}
  &\le |y_{k,0}P(\xi)| + |p| \dZ{y_{k,0}\xi^2}
       + |q| \dZ{y_{k,0}\xi} \cr
  &\le c_3 \big( Y_k H(\alpha) |\xi-\alpha| + Y_k^{-1} H(\alpha)
           \big) \cr}
$$
with a constant $c_3>0$ depending only on $\xi$.  Choosing $k$ to
be the smallest positive integer for which
$$
H(\alpha) \le {c_1\over 2c_3} Y_k^{1-\delta},
$$
and using (9.1) this implies
$$
|\xi-\alpha| \ge {c_1\over 2c_3} Y_k^{-1-\delta}H(\alpha)^{-1}.
$$
The conclusion follows since the above choice of $k$ implies $Y_k
\ll H(\alpha)^{\gamma/(1-\delta)}$. \cqfd

\medskip
Numerical experiments done on Fibonacci continued fractions
suggest that the values of $\dZ{y_{k,0}\xi^3}$ are more or less
uniformly distributed in the interval $(0,1/2)$.  If true, this
would imply a lower bound of the form $\dZ{y_{k,0}\xi^3}\gg 1/k
\gg (\log\log Y_k)^{-1}$.  To be safer, we may hope for a lower
bound of the type $\dZ{y_{k,0}\xi^3}\gg Y_k^{-\delta}$ for any
$\delta>0$. Then, for any $\epsilon > 0$, the above proposition
gives a measure of approximation by algebraic integers $\alpha$ of
degree $\le 3$ of the form  $|\xi-\alpha| \ge c
H(\alpha)^{-\gamma^2-\epsilon}$ with a constant $c$ depending on
$\xi$ and $\epsilon$.  Here we content ourselves with the
following much weaker lower bound which, in the same way, implies
Theorem 1.5.

\proclaim Proposition 9.2. In the above notation, we have
$\dZ{y_{k,0}\xi^3} \gg Y_k^{-1/\gamma^3}$ for all sufficiently
large values of $k$ (to avoid those indices with $y_{k,0}=0$).
\par

\proof  Fix an integer $k\ge 4$ and denote by $y_{k,3}$ the
nearest integer to $y_{k,2}\xi$.  Define
$$
\tuy_k=(y_{k,1},y_{k,2},y_{k,3}),
  \quad
d_{k-2}=\det(\uy_{k-2},\uy_{k-1},\uy_k)
  \and
\td_{k-2}=\det(\uy_{k-2},\uy_{k-1},\tuy_k,).
$$
According to Theorem 5.1, $d_{k-2}$ belongs to a finite set of
non-zero integers. Since $\xi$ is irrational and since $\td_{k-2}$
is an integer, we deduce that
$$
1 \ll |d_{k-2}\xi-\td_{k-2}|.
$$
On the other hand, using Proposition 8.1, we find
$$
\eqalign{ |d_{k-2}\xi-\td_{k-2}|
  &= |\det(\uy_{k-2}, \uy_{k-1}, \xi\uy_k-\tuy_k)| \cr
  &\ll \|\xi\uy_k-\tuy_k\| H(Q_{k-2}) \cr
  &\ll \big(\dZ{y_{k,2}\xi}+Y_k^{-1}\big) Y_{k-3}. \cr}
$$
Combining these two estimates, we obtain $\dZ{y_{k,2}\xi} \gg
Y_{k-3}^{-1} \gg Y_k^{-1/\gamma^3}$ for all sufficiently large
$k$.  The conclusion follows since $|y_{k,0}\xi^3-y_{k,2}\xi| \ll
Y_k^{-1}$. \cqfd

%
%

\section{10.  A partial converse}

We conclude this paper with the following partial converse to
Proposition 9.1, where the notation is as in \S9.

\proclaim Proposition 10.1.  Let $\ell$ be the least common
multiple of all integers of the form $d_k = \det(\uy_k, \uy_{k+1},
\uy_{k+2})$ with $k\ge 1$, and assume that there exist infinitely
many indices $k\ge 1$ satisfying
$$
\dZ{{y_{k,0}\xi^3\over \ell}} \le c_1 Y_k^{-\delta} \eqno(10.1)
$$
for some real numbers $c_1>0$ and $\delta>0$.  Then there exist a
constant $c_2>0$ and infinitely many algebraic integers
$\alpha\in\bC$ of degree $\le 3$ for which
$$
|\xi-\alpha| \le c_2 H(\alpha)^{-\theta}
$$
where $\theta=(\gamma^2-\delta)/(1-\delta)$.
\par

\medskip
Our proof below goes by making explicit the arguments of Davenport
and Schmidt in \S2 of [\DSb].  To make the connection clearer for
the interested reader, denote by $\cC_k$ the convex body of
$\bR^3$ defined by
$$
|x_0|\le Y_k^\gamma,
  \quad
|x_0\xi-x_1|\le Y_k^{-1}
  \and
|x_0\xi^2-x_2|\le Y_k^{-1}.
$$
Then there exists a constant $c>0$ which is independent of $k$
such that $\uy_{k+1}\in c\,\cC_k$, $\uy_{k}\in c\,\cC_k$ and
$\uy_{k-1}\in cY_k^{1/\gamma^2}\cC_k$ for $k\ge 2$.  Since the
volume of $\cC_k$ is $8Y_k^{-1/\gamma^2}$ and since $\uy_{k+1}$,
$\uy_{k}$ and $\uy_{k-1}$ are linearly independent over $\bQ$,
this implies that the successive minima of $\cC_k$ behave like
$1$, $1$ and $Y_k^{1/\gamma^2}$ and that they are essentially
realized by these three points. Then, the successive minima of the
polar convex body $\cC^*_k$ defined by
$$
|x_0+x_1\xi+x_2\xi^2| \le Y_k^{-\gamma},
  \quad
|x_1|\le Y_k
  \and
|x_2|\le Y_k
$$
behave like $Y_k^{-1/\gamma^2}$, $1$ and $1$ and, by identifying a
point $(x_0,x_1,x_2)\in\bZ^3$ with the polynomial $x_0+x_1T+x_2T^2
\in \bZ[T]$, these minima are essentially realized by the
polynomials denoted $B$, $A$ and $C$ in the proof below.

\proof Fix an index $k\ge 2$ for which (10.1) is satisfied.
Define $\uT=(1,T,T^2)$ and consider the three polynomials
$$
A=\det(\uT,\uy_{k-1},\uy_{k+1}),
  \quad
B=\det(\uT,\uy_{k},\uy_{k+1})
  \and
C=\det(\uT,\uy_{k-1},\uy_{k}).
$$
We find
$$
d_{k-1}
  = y_{k,0} A(T) + y_{k-1,0} B(T) + y_{k+1,0} C(T).
$$
Thus, if we put $m=\ell/d_{k-1}$, then, for a suitable choice of
signs, the polynomial
$$
P(T) = T^3 - \xi^3
       \pm m \dZ{{y_{k,0}\xi^3\over \ell}}A(T)
       \pm m \dZ{{y_{k-1,0}\xi^3\over \ell}}B(T)
       \pm m \dZ{{y_{k+1,0}\xi^3\over \ell}}C(T)
$$
has integral coefficients and is monic of degree three.   We claim
that it satisfies also
$$
|P'(\xi)|\sim H(P)
  \and
|P(\xi)| \ll H(P)^{-\gamma/(1-\delta)}
$$
and that $H(P)$ tends to infinity with $k$. If we take this for
granted, then the root $\alpha$ of $P$ which is closest to $\xi$
is an algebraic integer of degree $\le 3$ with
$$
|\xi-\alpha|
  \ll {|P(\xi)|\over |P'(\xi)|}
  \ll H(P)^{-\theta}
  \ll H(\alpha)^{-\theta}
$$
and the proposition is proved.

To establish the claim, we first note that, by (10.1) and
Proposition 9.2, we have
$$
Y_k^{-\delta}
  \gg \dZ{{y_{k,0}\xi^3\over \ell}}
  \ge {1\over \ell} \dZ{y_{k,0}\xi^3}
  \gg Y_k^{-1/\gamma^3},
\eqno(10.2)
$$
and so $\delta\le \gamma^{-3}$ since $k$ may be arbitrarily large.
We also observe that, in the notation of Proposition 8.1, we have
$B=Q_{k}$ and $C=Q_{k-1}$.  Thus, with the appropriate choice of
sign, we get
$$
H\left(P \pm m \dZ{{y_{k,0}\xi^3\over \ell}}A\right)
  \ll 1+H(B)+H(C)
  \ll Y_{k-1}
\eqno(10.3)
$$
and
$$
|P(\xi)| \ll |A(\xi)|+|B(\xi)|+|C(\xi)|
   \ll |A(\xi)|+Y_{k+1}^{-1}.
\eqno(10.4)
$$
Moreover, by Corollary 5.3, the point $d_{k-1}\uy_{k+2} - d_k
\uy_{k-1}$ is an integral multiple of $\uy_{k+1}$ and thus we find
$$
d_{k-1}Q_{k+1} + d_k A
  = \det(\uT, \uy_{k+1}, d_{k-1}\uy_{k+2} - d_k \uy_{k-1})
  = 0.
$$
Since $d_{k-1}$ and $d_k$ are non-zero integers of bounded
absolute value, this implies, by Proposition 8.1, that
$$
H(A) \sim |A'(\xi)| \sim Y_k
  \and
|A(\xi)| \ll Y_{k+3}^{-1}.
  \eqno(10.5)
$$
Combining (10.2), (10.3) and (10.5) we deduce
$$
H(P) \sim \dZ{{y_{k,0}\xi^3\over \ell}} Y_k \sim |P'(\xi)|.
$$
In particular this gives $H(P) \ll Y_k^{1-\delta}$ and so, using
(10.4) and (10.5), we get
$$
|P(\xi)|
  \ll Y_k^{-\gamma}
  \ll H(P)^{-\gamma/(1-\delta)}
$$
which completes the proof of the claim.  \cqfd

\proclaim Corollary 10.2. Let the hypotheses be as in Theorem 6.2
and assume moreover that $\trace(JAB)=\pm1$.  Then the following
conditions are equivalent:
\item{(i)} There exists $\theta>\gamma^2$ such that the inequality
$|\xi-\alpha| \le H(\alpha)^{-\theta}$ has infinitely many
solutions in algebraic integers $\alpha\in\bC$ of degree $\le 3$
over $\bQ$.
\item{(ii)} There exists $\delta>0$ such that the inequality
$\{y_{k,0}\xi^3\} \le Y_k^{-\delta}$ holds for infinitely many
indices $k$.
\par

\proof The fact that (i) implies (ii) follows from Proposition
9.1.  The reverse implication follows from Proposition 10.1
granted that, for the sequence $(\uy_k)_{k\ge -1}$ constructed by
Lemma 6.1, we have $|\det(\uy_k, \uy_{k+1}, \uy_{k+2})|=1$ for all
$k\ge -1$.  \cqfd

\medskip
In particular, the corollary applies to any Fibonacci continued
fraction $\xi_{a,b}$ for which $|a-b|=1$.

\medskip
\noindent{\bf Acknowledgments.} The author thanks Drew Vandeth for
the reference [\ADQZ] and for recognizing the continued fraction
expansion of Corollary 6.3 as a Fibonacci word.

%
%
\penalty-1000
\section{References}

\medskip
\item{[\ADQZ]} J.-P.~Allouche, J.~L.~Davison, M.~Queff\'elec,
L.~Q.~Zamboni, Transcendence of Sturmian or morphic continued
fractions, {\it J.~Number Theory} {\bf 91} (2001), 39--66.

\item{[\BT]} Y.~Bugeaud, O.~Teuli\'e, Approximation d'un nombre
r\'eel par des nombres alg\'ebriques de degr\'e donn\'e, {\it Acta
Arith.\ }{\bf 93} (2000), 77--86.

\item{[\Ca]} J.~W.~S.~Cassels, {\it An introduction to Diophantine
approximation}, Cambridge Univ.\ Press, 1957.

\item{[\DSa]} H.~Davenport, W.~M.~Schmidt,  Approximation to real
numbers by quadratic irrationals, {\it Acta Arith.\ }{\bf 13}
(1967), 169-176.

\item{[\DSb]} H.~Davenport, W.~M.~Schmidt,
Approximation to real numbers by algebraic integers, {\it Acta
Arith.\ }{\bf 15} (1969), 393--416.

\item{[\LR]} M.~Laurent, D.~Roy, Criteria of algebraic independence
with multiplicities and interpolation determinants, {\sl Trans.\
Amer.\ Math.\ Soc.} {\bf 351} (1999), 1845--1870.

\item{[\Lo]} M.~Lothaire, {\it Combinatorics on words}, Encyclopedia
of mathematics and its applications, vol.~17, Addison-Wesley Pub.\
Co., 1983.

\item{[\Ro]} D.~Roy, Approximation simultan\'ee d'un nombre et de son
carr\'e, {\it C.\ R.\ Acad.\ Sc.\ Paris, S\'erie I}, {\bf 336}
(2003), 1--6, arXiv:math.NT/0210395.

\item{[\Sc]} W.~M.~Schmidt, {\it Diophantine approximation}, Lecture
Note in Math., vol.~785, Sprin\-ger-Verlag, 1980.

\item{[\Wi]} E.~Wirsing, Approximation mit algebraischen Zahlen
beschr\"ankten Grades, {\it J.\ reine angew.\ Math.\ }{\bf 206}
(1961), 67-77.

 \vskip 2truecm plus .5truecm minus .5truecm

\hbox{ \vtop{
 \hbox{Damien ROY}
 \hbox{D\'epartement de Math\'ematiques et de Statistiques}
 \hbox{Universit\'e d'Ottawa}
 \hbox{585 King Edward}
 \hbox{Ottawa, Ontario K1N 6N5, Canada}
 \hbox{{\it E-mail:} droy@uottawa.ca}
 \hbox{{
http{$:$}//aix1.uottawa.ca/${\scriptstyle \sim}$droy/}} }}

\bye